\documentclass[times,10pt,twocolumn]{article} 
\usepackage{latex8}
\usepackage{times}
\usepackage{graphics}
\usepackage{latexsym}

\begin{document}
\newcommand{\tab}[0]{\hspace{.1in}}
\newtheorem{thm}{\noindent Theorem}
\newtheorem{lem}{\noindent Lemma}
\newenvironment{proof}{\noindent {\bf Proof }}{\vspace{.1in}}



\vspace{-2in}
\title{The Randomness Recycler:  A New Technique for Perfect Sampling}
\author{James Allen Fill\thanks{Research 
supported by NSF grant DMS-9803780
and by the Acheson J. Duncan Fund for the Advancement of Research in
Statistics.}
\\Department of Mathematical Sciences
\\The Johns Hopkins University\\
jimfill@jhu.edu\\
\and
Mark Huber\thanks{Supported by NSF postdoctoral fellowship 99-71064}
\\Department of Statistics
\\Stanford University\\
mhuber@orie.cornell.edu\\
} 
\maketitle
\thispagestyle{empty}

\begin{abstract}
For many probability distributions of interest, it is quite
difficult to 
obtain samples efficiently.  Often,
Markov chains are employed
to obtain approximately random samples from these
distributions.  The primary
drawback to traditional Markov chain methods is that the mixing
time of the chain is usually unknown, which makes 
 it impossible to
determine how close the output samples
are to having the 
target distribution.
Here we present a new protocol, the randomness recycler (RR), that  
overcomes this difficulty.  Unlike classical 
Markov chain approaches, an RR-based algorithm creates samples
drawn exactly from the desired distribution.  Other perfect
sampling methods such as coupling from the past use
existing Markov chains, but RR
does not use the traditional Markov chain at all.
While by no means universally useful, RR does
apply to a wide variety of problems.  In restricted
instances of certain problems, it gives the first expected
linear time algorithms for generating samples.  Here we
apply RR to self-organizing lists, the Ising model, random
independent sets, random colorings, and the random cluster model.
\end{abstract}

\section{Introduction}
\label{SEC:I}
\pagenumbering{arabic}

The Markov chain Monte Carlo (MCMC) approach to generating samples
has enjoyed enormous success since its introduction, but in
certain cases it is possible to do better.  
The ``randomness recycler'' technique we introduce here
(and whose name is explained in Section~\ref{SEC:WIS}) works
for a variety of problems without employing the traditional Markov 
chain.  Our approach is faster in many cases,
generating in particular the first algorithms that have expected running time
linear in the size of the problem, under certain restrictions.

In classical MCMC approaches, small random changes are made 
in the observation 
until the observation has nearly the stationary distribution
of the chain.  The Metropolis~\cite{metropolisrrtt53} and heat 
bath algorithms utilize the idea of reversibility to design chains with a
stationary distribution matching the desired distribution.
Unfortunately, this standard Markov chain approach does have problems.

The samples generated by MCMC will not be drawn exactly
from the stationary distribution, but only approximately.
Moreover, they will not be close to the stationary distribution
until a number of steps larger than the mixing time of the chain
have been taken.  Often the mixing
time is unknown, and so the quality of the sample is suspect.

Recently, Propp and Wilson have shown how to avoid these problems
using techniques such as coupling from the past (CFTP)~\cite{proppw96}.
For some chains, CFTP provides a procedure that allows perfect samples 
to be drawn from the stationary distribution of the chain, without
knowledge of the mixing time.  However, CFTP and related 
approaches have drawbacks of their own.  
These algorithms are
noninterruptible, which means that the user must commit to 
running such an algorithm
for its entire (random) running time even though that time is not
known in advance.  Failure to do so can introduce bias
into the sample.
Other algorithms, such as FMMR~\cite{fill99}, are interruptible
(when time is measured in Markov chain steps), but
require storage and subsequent rereading of random bits used by the algorithm.  
The method we will present is both interruptible and ``read-once,'' with no
storage of random bits needed.

In addition,
algorithms like CFTP and FMMR require an underlying Markov chain,
and can never be faster than the mixing time of this 
underlying chain.  Often these chains
make changes to parts of the state where the state has already been 
suitably randomized.  This leads to wasted effort when running the algorithm
that often adds a log factor to the running time of the algorithm.

The randomness recycler (RR) is not like any of these perfect sampling
algorithms.  In fact, the RR approach 
abandons the traditional Markov chain 
entirely.  This is what allows the algorithm in several cases to
reach an expected running time that is linear, the first for several
problems of interest.  The RR technique gives interruptible, 
read-once perfect samples.

In the next section we illustrate the randomness recycler
for the problem of finding random independent sets of
a graph.  After this example we present in Section~\ref{SEC:RR} the general
randomness recycler procedure
and present a (partial) proof of correctness.  
In Section~\ref{SEC:Appl} we present other applications
and in Section~\ref{SEC:Conc} we review the
results of applying our new approach to several different problems.

\section{Weighted Independent Sets}
\label{SEC:WIS}
We begin by showing how the randomness recycler 
technique applies to the problem
of generating a random independent set of a graph.  This
will illustrate the key features of RR and lay the groundwork for
the more general procedure described in the next section.

Recall that an \emph{independent set} of a graph is a subset of vertices no two
of which share an edge.  We will represent an independent set as a 
coloring of the vertices from $\{0,1\}$, denoted generically by~$x$.
Set
$x(v) = 1$ if $v$ is in the independent set, and
$x(v) = 0$ if $v$ is not in the independent set.  
This implies that $\sum_v x(v)$ is
 the size of the independent set.

We wish to sample from the distribution
$$\pi(x) = \frac{\lambda^{\sum_v x(v)}}{Z_\lambda},$$
where $\lambda$ (called the \emph{fugacity}) is a parameter of the problem, and
$Z_\lambda$ is the normalizing constant needed to make~$\pi$ a probability
distribution.

This distribution is known as the hard core gas model 
in statistical physics,
and also has applications in stochastic loss networks~\cite{kelly91}.
When $\lambda$ is large the sample tends to be a large independent set,
and if $\lambda > 25 / \Delta$ where 
$\Delta$ is the maximum degree of
the graph, it is known that generating samples from this distribution 
cannot be done in polynomial time unless $NP = RP$~\cite{dyerfj98}.

We will show that for $\lambda < 4 / (3 \Delta - 4)$  
the randomness
recycler approach gives an algorithm with expected running time 
linear in the size of the graph, the first such result for this problem.

The RR approach is to start not with the entire graph, but rather with a small
graph where we can easily find an independent set from this distribution.
For example, if a graph has only a single vertex, finding an independent set is
easy.  Starting from a single vertex, we attempt to add vertices to the graph,
building up until we are back at our original problem.  Sometimes we fail in our
attempt to build up the graph, and indeed will then also need to remove vertices
that we had previously added.  The set
$V_t$ will comprise those vertices we have built up by the end of time step~$t$.
After step~$t$, the vector $x_t$ will hold an independent set which has the
correct distribution over (the subgraph induced by) the vertices in $V_t$.

\begin{center}
\vspace{.1in}
\begin{tabular}{|l|}
\hline
{\bf Randomness Recycler for}
{\bf Independent Sets} \\
\\ 
{\bf Set} $V_0 \leftarrow \emptyset$, $x_0 \leftarrow 0$, $t \leftarrow
0$ \\
{\bf Repeat} \\
\tab {\bf Set} $x_{t+1} \leftarrow x_t$ \\
\tab {\bf Choose} any $v \in V \setminus V_t$ \\
\tab {\bf Set} $V_{t+1} \leftarrow V_t \cup \{v\}$ \\
\tab {\bf Draw} $U$ uniformly at random from $[0,1]$ \\
\tab {\bf If} $U \leq 1 / (1 + \lambda)$ \\
\tab \tab {\bf Let} $x_{t+1}(v) \leftarrow 0$ \\
\tab {\bf Else} \\
\tab \tab {\bf Let} $x_{t+1}(v) \leftarrow 1$ \\
\tab \tab {\bf If} a neighbor $w$ of $v$ has 
		$x_{t+1}(w) = 1$ \\
\tab \tab \tab {\bf Let} $w$ be the lowest-numbered 
		such neighbor \\
\tab \tab \tab {\bf Set} $x_{t+1}(w) \leftarrow 0$, 
   $x_{t+1}(v) \leftarrow 0$ \\
\tab \tab \tab {\bf Remove} from $V_{t+1}$ the vertices $v$ and $w$,\\
\tab \tab \tab \tab  all neighbors of~$w$, and all neighbors of~$v$ \\
\tab \tab \tab \tab  with numbers less than that of~$w$ \\
\tab {\bf Set} $t \leftarrow t+1$ \\
{\bf Until} $V_t = V$ \\
\hline
\end{tabular}
\end{center}

(In advance of running the algorithm, choose and fix a numbering of the vertices.)
The algorithm proceeds inductively as follows.
At the outset of step $t + 1$,
we begin with
an independent set $x_t$ of $V_t$ chosen with the correct probability.
Then we choose a vertex~$v$ not in 
$V_t$ to attempt to add.  This vertex may be chosen in any fashion
desired (randomly, or according to some fixed order, but not depending on the
independent set~$x_t$).  Because the desired probability of choosing an
independent set~$x$ is proportional to
$\lambda^{\sum x(v)}$, putting $x_{t+1}(v) = 1$ has $\lambda$ times the
weight of putting $x_{t+1}(v) = 0$.  Therefore we select $x_{t+1}(v) = 1$
with probability $\lambda / (1 + \lambda)$ and $x_{t+1}(v) = 0$ with
probability $1 / (1 + \lambda)$ (these are the heat
bath probabilities).

Unfortunately, the vector~$x_{t + 1}$ resulting from 
this selection may fail to correspond to
an independent set.  At line 11 of the pseudocode, we check whether
some neighbor of $v$ was already colored 1 (in
the independent set).  Note that we cannot
simply remove $v$.  Prior to the step, we knew that $x_t$ was an independent
set of $V_t$.  If we observe that $x_t(w) = 1$ for some
lowest-numbered neighbor~$w$ of
$v$, then $x_t$ is an independent
set on $V_t$ conditioned on this knowledge.

Our solution is this:\ \ In line 14 we 
``undo'' the knowledge gained by removing from $V_{t+1}$ the vertices $v$ and $w$,
all the neighbors of~$w$, and all the neighbors of~$v$ with number less than that
of~$w$.  On the remaining vertices of $V_{t+1}$, 
$x_{t + 1}$ will continue to be an independent set from the correct distribution.
We will say that an RR step of this type preserves the correct distribution.

Note that although $V_{t+1}$ is made smaller than $V_t$ in the case of a
conflict, we are able to salvage most of the vertices in $V_t$.  In other words,
we ``recycle'' the randomness built up in all of the vertices except $v$ and $w$
and some neighbors. This is where our approach gets its name, 
and ``recycling''
is the key new feature that enables us to 
contruct similar practicable algorithms
for a wide variety of problems.

We repeat until $V_t = V$.  Because each step preserves the correct
distribution, we know that $x_t$ will have the correct
distribution $\pi$ at the end.  This is proved formally in the
next section; here we concentrate on bounding the running time of our
procedure.

\begin{thm}
\label{THM:WIS}
If $\lambda < 1 / (2\Delta - 1)$, then the expected running time of the above
randomness recycling procedure for random independent sets is $O(n)$.
\end{thm}

A more careful statement of Theorem~\ref{THM:WIS} 
is given following the proof.
\smallskip

\begin{proof}
We will show that for $\lambda$ this small, on
average $|V_t|$ increases at each step.  If $U \leq 1 / (1 + \lambda)$,
then the size of $|V_t|$ goes up by~$1$, but if $U > 1 / (1 + \lambda)$,
then the size of $|V_t|$ may decrease by at most $2 \Delta - 1$ [removing
$v$ (not previously included), $w$, and some neighbors].  Hence
\begin{eqnarray*}
E \bigl[ \bigl. |V_{t+1}| - |V_t|\,\bigr|\, V_t, x_t \bigr]
 & \geq & \frac{1}{1 + \lambda} (1)
  - \frac{\lambda}{1 + \lambda} (2 \Delta - 1) \\
& = & \frac{1}{1 + \lambda} [1 - (2 \Delta - 1)\lambda],
\end{eqnarray*}
which is positive precisely when $\lambda < 1 / (2 \Delta - 1)$.  Given an
increase of $|V_t|$ on average at each step, standard martingale stopping
theorems (see, e.g.,\ \cite{port94}) show that after $O(n)$ 
expected time the
value of $|V_t|$ will be $n$, at which point $V_t = V$ and the algorithm
terminates. $\Box$
\end{proof}

More carefully, if
$$
\frac{1}{1 + \lambda} [1 - (2 \Delta - 1) \lambda] \geq \gamma \in (0, 1),
$$
i.e.,\ if
$$
\lambda \leq 1 \left/ \left( \frac{2 \Delta}{1 - \gamma} - 1 \right) \right.,
$$
then the expected value of the running time~$T$ (measured by number of
iterations of the Repeat loop) satisfies
$$
E T \leq n / \gamma.
$$
Furthermore, a simple argument shows that the distribution of~$T$ has at worst
geometrically thick tails:
$$
P(T \geq 2 \mbox{$\frac{m}{\gamma}$} n) \leq 2^{-m},\ \ \ m = 1, 2, \ldots.
$$

Several tricks may be used to either to improve our method or to improve our
bounds on its performance.  The first two we present work by altering the
algorithm, and the third gives a better analysis. First we concentrate on making
sure that as few vertices as possible are removed in the rejection step.  Note
that we may assume that the graph is connected, since otherwise we simply work on
each connected component separately.  Therefore $V \setminus V_0$ is connected,
and by being slightly careful in how we choose $v \in V \setminus V_t$, we
can ensure that $V \setminus V_t$ remains connected at each step.  
In every step (except when $|V \setminus V_t| = 1$), the vertex $v$ is
not adjacent to $\Delta$ vertices in $V_t$, but only to at most $\Delta - 1$, so
fewer vertices are removed during rejection.  Since the vertices removed from 
$V_t$ in case of rejection are connected to $V \setminus V_t$, 
$V \setminus V_{t+1}$ will also be connected whether we accept or reject,
and no more than an extra constant amount of work is required at each step.

The second alteration concerns how we look for the neighbor of $v$ that
is colored 1 in case of rejection.  Instead of starting at the lowest
numbered neighbor and working our way up, we start at a random neighbor
and continue looking in cyclical order until we find our~$w$ colored~1;
and then the vertices that we remove are $w$ and its neighbors and $v$
and its neighbors encountered in the search prior to finding~$w$.
On average (and together with the first trick),
we need only look at $\leq \Delta / 2$ vertices in order to find~$w$.

Finally, in our analysis we kept track of a potential 
$\phi(V_t, x_t) = |V_t|$ and showed that~$\phi$ increases on average.
When we accept we sometimes add a vertex colored~1 to our set~$V_t$;
but when we reject, precisely one vertex colored~1 (namely, $w$) is removed.
This suggests that we modify~$\phi$
so that the acceptance and rejection phases both lead
(in the worst case)
to the 
same expected change in $\phi$.
We will consider
$$\phi(V_t, x_t) = |V_t| - \alpha \sum_{v} x_t(v)$$
and seek a suitable value of~$\alpha$.
The expected change in $\phi$ if no neighbor of $v$ is colored~1 is
$1 - \alpha \left[ \lambda / (1 + \lambda) \right]$.  If some neighbor
is colored~1, then the expected change is at least
$1 \left[1 / (1 + \lambda) \right] + 
  \left(\frac{\lambda}{1 + \lambda}\right) \left[ 
  - \left(\frac{3 \Delta - 2}{2} \right) + \alpha \right].$
[The term $(3 \Delta - 2) / 2$ is an upper bound on the expected decrease
in $|V_t|$, since (see above) on average we lose at most $\Delta / 2$ neighbors
of~$v$ and $\Delta - 1$ neighbors of~$w$.]

These two expressions may be made equal by setting 
$\alpha = 3 \Delta / 4$, and then the expected change in $\phi$ will
be positive when 
$$\lambda < \frac{4}{3 \Delta - 4}.$$

Note that $\phi$ at time 0 equals 0, and can never be more than $n$, 
and we have shown that $\phi$ is expected to increase by a fixed positive
amount at each step (when we formulate carefully as in the paragraph following
the proof of Theorem~\ref{THM:WIS}).  This fact together with standard martingale
stopping theorems can then be used to show that the expected time 
needed for $|V_t|$ to equal~$n$ is at most linear in $n$.

\subsection{Markov chain approaches}
Several Markov chains for this problem 
exist~\cite{lubyv99}~\cite{dyerg97a}, together with techniques for using
CFTP to obtain perfect samples~\cite{haggstromn98}~\cite{huber98}.  
These Markov chains are known to mix in time $O(n \ln n)$, and the
corresponding perfect sampling algorithms are known to run in time
$O(n \ln n)$, when $\lambda < 2 / (\Delta - 2)$,
which is a larger range of $\lambda$ than our method gives.  
However, when $\lambda$ is small enough, the $O(n)$ bound for our RR algorithm
is smaller.  It is hoped that with further refinement of the
rejection step, the range of $\lambda$ may be increased to where
it matches the Markov chain analysis.

\section{The Randomness Recycler}
\label{SEC:RR}
We now present a more general outline of the randomness recycler technique.
Many state spaces~$\Omega$ of interest are of the form
$\Omega \subseteq C^V$, where $C^V$ is the set of (proper or improper) 
colorings of a graph.  Our goal is to sample from $\Omega$ in expected
time linear in $|V|$.
We have already seen how the independent sets of
a graph may be encoded by coloring a vertex 1 if it is in the indpendent
set and 0 otherwise.  For another example, the set of permutations of $n$
elements is a subset of $\{1,\ldots,n\}^{\{1,\ldots,n\}}$.  Of course,
the size
of $\Omega$ may be as large as $|C|^{|V|}$, and this is in part what makes
generating samples from these distributions difficult.

\begin{center}
\vspace{.1in}
\begin{tabular}{|l|}
\hline
{\bf The Randomness Recycler (Outline)}\\
\\
{\bf Set} $V_0 \leftarrow \emptyset$, $X_0 \leftarrow \mbox{suitable~$x_0$}$, $t
\leftarrow 0$ \\
{\bf Repeat} \\
\tab {\bf Set} $X_{t+1} \leftarrow X_t$ \\
\tab {\bf Choose} $v \in V \setminus V_t$ \\
\tab {\bf Randomly} choose color $c$ for $v$ \\
\tab {\bf Compute} probability of accepting color $c$ for $v$ \\
\tab {\bf If} we accept \\
\tab \tab {\bf Set} $X_{t+1}(v) \leftarrow c$ \\
\tab \tab {\bf Set} $V_{t+1} \leftarrow V_t \cup \{v\}$ \\
\tab {\bf Else} \\
\tab \tab {\bf Set} $V_{t+1}$ and $X_{t+1}|_{V \setminus V_{t+1}}$ in a way that
  \\
\tab \tab \tab `undoes' the effect of rejection \\ 
\tab {\bf Set} $t \leftarrow t + 1$ \\
{\bf Until} $V_t = V$ \\
\hline
\end{tabular}
\end{center}

In an RR algorithm, a sample 
(i.e.,\ one draw from~$\pi$) is built up one vertex of $V$ at a time
until we include all of the vertices.
Let $V_t$ be the subset of vertices 
on which we have already built up a sample at time~$t$.  
On the vertices in
$V \setminus V_t$, the sample is fixed at some value, whereas on
$V_t$, the sample is random, and drawn exactly from the desired
distribution.  $V_t$ starts out empty, and at each step of the
algorithm we attempt to add a vertex to $V_t$.  Sometimes this is possible,
and sometimes it is not.  We continue in this fashion until $V_t = V$,
at which point we have a sample drawn exactly from the desired 
distribution.
Let $X_t$ denote the coloring of the graph~$V$ at time $t$.


The way in which we randomly choose $c$, compute the acceptance
probability, and set $V_{t+1}$ and $X_{t+1}|_{V \setminus V_{t+1}}$ in case of
rejection will all depend on the target distribution~$\pi$.  What differentiates
this algorithm from an elementary stepwise rejection approach
is our rejection step.  Rather than starting over when rejection is faced,
we keep as much of $V_t$ as possible,
``recycling'' the coloring on $V_{t+1}$.  

At each time step $t$ we keep track of the vertex set $V_t$ 
together with the colors that are fixed on $V \setminus V_t$.
The state $X_t$ is random over $V_t$ while on $V \setminus V_t$
it is deterministic. 
Let $X^*_t = (V_t, X_t|_{V \setminus V_t})$, and
for any possible value $x^* = (S, x|_{V \setminus S})$ of $X^*_t$,
let $\pi_{x^*}$ be~$\pi$ conditionally given that the colors of
$V \setminus S$ are as specified by $x|_{V \setminus S}$.

To achieve both the desired distribution and interruptibility,
we want $X_t$ to be random over $V_t$ independent of the history 
$X^*_{t'}$ for $t' < t$.  In other words we want the identity
\begin{equation}
\label{eqn:sp}
P(X_t = x_t | X^*_0 = x^*_0, \ldots, X^*_{t} = x^*_{t})
  = \pi_{x^*_t}(x_t),
\end{equation}
to hold.  Indeed, if it does, then letting~$T$
denote the first time that $V_T = V$, it follows easily that
$$P(X_T = x|T = t) = \pi(x).$$
Thus if~(\ref{eqn:sp}) is satisfied for all~$t$, then 
at termination time~$T$ the RR algorithm
returns a sample~$X_T$ that is distributed according to the desired
distribution, and we have the interruptibility property that~$T$
and~$X_T$ are independent random variables.

Since $V_0$ is empty, it is easy to begin with $X_0$ from
$\pi_{x^*_0}$.  
Let
$H_t := (X^*_0 = x^*_0, X^*_1 = x^*_1, \ldots, X^*_t = x^*_t)$
for notational convenience. 
We will say that 
step $t + 1$ \emph{preserves the correct distribution} if
\begin{eqnarray*}
& P(X_t = x_t | H_t) \equiv \pi_{X^*_t}(x_t) \\
&  \Downarrow \\
&  P(X_{t+1} = x_{t+1} | H_{t+1}) \equiv \pi_{x^*_{t+1}}(x_{t+1}).
\end{eqnarray*}

This requirement that RR preserve the correct distribution
is somewhat analogous to the design requirement that a Markov chain be
reversible. It gives us a straightforward approach to designing an
RR.  

Just as the heat bath approach gives a means for designing Markov
chains that are reversible, it also gives us a method for 
designing RR algorithms that preserve the correct distribution.
For a specified vertex $v \in V$ and coloring~$x$,
let $\pi_v(\cdot; x)$ denote the conditional probability distribution
of~$X(v)$ given that $X|_{V \setminus \{v\}} = x|_{V \setminus \{v\}}$
when~$X$ has the stationary distribution~$\pi$. 
From current state~$x$, the heat bath (or Gibbs sampler)
Markov chain approach is to choose~$v$ uniformly at random and then
choose a new color for $v$ distributed according to 
$\pi_v(\cdot; x)$.

In heat bath RR, the vertex $v$ is 
chosen
any way the user desires
from $V \setminus V_t$, and then a new color is picked according
to $\pi_v(\cdot; x)$.  However, this color is not always accepted.  We
compute the acceptance probability as follows, with the goal being
to preserve the correct distribution.  According to Theorem~\ref{THM:correct}
below, this goal is indeed met.

Given values $x^*_t$, $x_t$, $x^*_{t+1}$, and $x_{t+1}$ that correspond to
a possible acceptance step in which vertex~$v$ is added to the growing vertex
set, define $\rho(x^*_t, x_t, x^*_{t+1}, x_{t+1})$ to be the ratio
$$
\rho(x^*_t, x_t, x^*_{t+1}, x_{t+1})
 := \frac{\pi_{x^*_{t+1}}(x_{t+1})}{\pi_v(x_{t+1}(v); x_t) \pi_{x^*_t}(x_t)}.
$$
Also define
$$
M(x^*_t, x^*_{t+1}) := \max_{x_t, x_{t+1}} \rho(x^*_t, x_t, x^*_{t+1}, x_{t+1}).
$$
Then the probability that we accept a possible transition from $(x^*_t, x_t)$ to
$(x^*_{t+1}, x_{t+1})$ is taken to be $\rho(x^*_t, x_t, x^*_{t+1}, x_{t+1}) /
M(x^*_t, x^*_{t+1})$. 

We do not have to use the heat bath probabilities.
It is also valid to use the same acceptance probability, with the distributions
$\pi_v(\cdot; x)$ replaced by arbitrary distributions $p_v(\cdot; x)$, when
the distribution $p_v(\cdot; x)$ is used to color a selected~$v$ when at a
configuration~$x$.

While these acceptance probabilities may appear daunting, for many problems
they simplify considerably.  For instance, in the independent
set case, suppose first that~$v$ has no neighbor colored 1.  Then the
heat bath probabilities are $1 / (1 + \lambda)$ for color~0
and $\lambda / (1 + \lambda)$ for color~1.  The acceptance probability
in this first case will always be 1.

If instead some neighbor of~$v$ is colored 1, then heat bath assigns
probability~$1$ to the color~$1$.  The acceptance probability, however, works out
to $1 / (1 + \lambda)$.  Careful examination of the independent set algorithm
in Section~\ref{SEC:WIS} shows that this is exactly how the color for $v$ is
chosen, with the same acceptance probabilities.

To show that the heat bath 
randomness recycler approach actually works (in general),
we need to show that every step preserves the correct distribution.
We will first consider acceptance steps, for which
the following lemma gives a sufficient condition.
\begin{lem}
\label{lem:sim}
Given possible values
$x^*_t$, $x^*_{t+1}$, and $x_{t+1}$ of
$X^*_t$, $X^*_{t+1}$, and $X_{t+1}$
corresponding to an acceptance step,
suppose that only one value $x_t$ of
$X_t$ has positive probability.  
If the bivariate process $(X^*_t, X_t)_{t \geq 0}$ evolves Markovianly
and if
for all such $x^*_t$, $x^*_{t+1}$, and $x_{t+1}$
and the single $x_t$ they determine we have
\begin{eqnarray*}
\lefteqn{\hspace{-1in}P(X^*_{t+1} = x^*_{t+1}, X_{t+1} = x_{t+1} | 
	X^*_t = x^*_t, X_t = x_t) \pi_{x^*_t}(x_t)} \\
& = & \pi_{x^*_{t+1}}(x_{t+1}) C,
\end{eqnarray*}
where $C$ does not depend on $x_t$ or $x_{t+1}$, then step $t + 1$ preserves the
correct distribution.
\end{lem}

\begin{proof}
Let $C_1 := 1 / P(X^*_{t+1} = x^*_{t+1}| H_t)$, and suppose that
$P(X_t = x_t| H_t) = \pi_{x^*_t}(x_t)$.  Let $E$ be the event
that $X^*_{t+1} = x^*_{t+1}$ and $X_{t+1} = x_{t+1}$.  Then
\begin{eqnarray*}
\lefteqn{P(X_{t+1} = x_{t+1} | H_{t+1} )} \\
& = &
	C_1 P(X^*_{t+1} = x^*_{t+1}, X_{t+1} = x_{t+1} | H_t)  \\
& = & C_1 P(E \cap
 \{X_t = x_t\} |H_t) \\
& = & C_1 P(X_t = x_t | H_t) P(E
	| H_t \cap \{X_t = x_t\}) \\
& = & C_1 \pi_{x^*_t}(x_t) P(E
 | X^*_t = x^*_t, X_t = x_t) \\
& = & C_1 C \pi_{x^*_{t+1}}(x_{t+1}),
\end{eqnarray*}
where the last step is exactly our assumption.

Note that neither~$C_1$ nor~$C$ depends on~$x_{t+1}$.
Hence, summing over~$x_{t + 1}$, 
$$1 = C_1 C \sum_{x_{t+1}} \pi_{x^*_{t+1}}(x_{t+1}) = C_1 C.$$
This completes the proof. $\Box$
\end{proof}

\begin{thm}
\label{THM:correct}
The heat bath RR and arbitrary RR acceptance steps 
preserve the correct distribution.
\end{thm}

\begin{proof}
The acceptance probabilities were chosen
precisely to match the requirements of Lemma~\ref{lem:sim}.
For instance, with heat bath RR, the left side of the equation in
Lemma~\ref{lem:sim} equals
\begin{eqnarray*}
\lefteqn{\hspace{-1.3in}\pi_v(x_{t+1}(v); x_t) \times
\frac{\pi_{x^*_{t+1}}(x_{t+1})}{\pi_v(x_{t+1}(v); x_t) \pi_{x^*_t}(x_t) M(x^*_t,
x^*_{t+1})}} \\
& & \hspace{-.7in}{} \times \pi_{x^*_t}(x_t),
\end{eqnarray*}
which reduces to the right side of the equation with $C = 1 / M(x^*_t,
x^*_{t+1})$.  The calculation for arbitrary RR is entirely similar.~$\Box$
\end{proof}

Now we turn our attention to rejection steps.
In designing an RR algorithm, it is our experience that proper handling
of rejection steps to ensure preservation of the correct distribution
is more difficult and problem-specific 
to arrange than is proper handling of acceptance steps.
But here are some broad guiding comments. 

Determination of the acceptance probability at step $t + 1$
will reveal knowledge about the colors of some subset, call it~$D_t$, of~$V_t$.
If we reject, we then set $V_{t+1}$ to be $V_t \setminus D_t$.
This insures that when we reject,  we do not bias the sample.
That is, by removing~$D_t$ from~$V_t$, we remove all traces of our knowledge
gained, and as a result the remaining sample is drawn exactly from
$\pi_{V_{t+1}}$.

In the case of the independent sets, the set~$D_t$ consists of
precisely those vertices prescribed to be removed by the algorithm:\ 
$w$ and all its neighbors, and neighbors of~$v$ with numbers lower than that
of~$w$.  [Indeed, all of these vertices are colored~0 at time~$t$, except for
vertex~$w$, which is colored~1.] 
It is not hard to check rigorously in this case that rejection steps
also preserve the correct distribution, but we omit the details.

\section{Applications}
\label{SEC:Appl}
This section applies the
randomness recycler approach to several different problems of
interest.  For some of these models we have theoretical bounds
on the running time, while for others we have only experimental results. 
 
\paragraph{The Ising and Potts models}
In the Ising model, vertices in a graph $(V,E)$ 
are colored from the set $\{-1,1\}$.
The distribution~$\pi$ from which we wish to sample is defined by
$$\pi(x) := \frac{\exp\{-\beta J H(x)\}}{Z_\beta},$$
where
$$H(x) := - \sum_{\{v_1,v_2\} \in E} x(v_1) x(v_2)$$
is known as the energy of the coloring,
$\beta$ is (proportional to) a postive parameter known as inverse temperature,
and~$J$ is~$1$ in the ferromagnetic model and~$-1$ in the antiferromagnetic
model.  Generating approximate samples may
be done in (nonlinear)
polynomial
time in the ferromagnetic case using Markov chain techniques
of Jerrum and Sinclair~\cite{jerrums93}~\cite{randallw98}.  

The RR approach has provably linear expected running time for
both the ferromagnetic and antiferromagnetic models when 
$\beta$ is small (i.e., the temperature is high).  The set $D_t$ to 
be removed from $V_t$ in case of rejection is just the set of neighbors of
the vertex~$v$ that we tried to add.  Omitting details and proofs, we
simply state the running time bound in the following theorem.
\begin{thm}
\label{THM:Ising}
Let $\Delta$ be the maximum degree of the graph.
If
$$e^{\beta} < \left( 1 + \frac{1}{\Delta} \right)^{1/\Delta},$$
then the expected running time of the heat bath RR procedure
for the Ising model is $O(n)$.
\end{thm}

Comments like those following the proof of Theorem~\ref{THM:WIS}
apply here, where now the expected increase
$\frac{1}{1 + \lambda} [1 - (2 \Delta - 1) \lambda]$
in~$|V_t|$ becomes $(\Delta + 1) e^{- \beta \Delta} - \Delta$. 

Although not needed for the theorem, in practice 
it helps to introduce a third color 0 to supplement
$\{-1,1\}$.  Notice that no edge with an endpoint colored~0 
contributes to $H$.  At the completion of step~$t$, every vertex in 
$V \setminus V_t$
which is surrounded entirely by vertices in $V \setminus V_t$ may be recolored 
0 since this action does not affect the vertices in $V_t$ at all.

The Potts model differs from the Ising model in that more than two colors are
used, but the energy depends (in a natural way) only on whether edges are colored
concordantly or discordantly, and the running time Theorem~\ref{THM:Ising}
remains valid verbatim.

\paragraph{The Random Cluster Model}
The random cluster model is an extension of the Potts model to noninteger
numbers of colors~\cite{fortuink72}; this is discussed further below.  
Unlike our previous examples, which colored vertices, the
random cluster model colors edges of a given graph $G = (V, E)$
with colors from $\{0,1\}$.  If $A$ is the set
of edges colored $1$, then the distribution is
$$\pi(A) := p^{|A|}(1 - p)^{|E \setminus A|} q^{c(A)} / Z_{p, q},
\ \ \ A \subseteq E,$$
where $p \in [0,1]$; $q > 0$ is not necessarily an integer, and we shall assume $q
> 1$; $c(A)$ is the number of connected components in the graph $(V, A)$; and
$Z_{p, q}$ is a normalizing constant.

The RR approach is as follows.  We represent a set $A \subseteq E$
by a binary vector~$x$, by setting $x(e) = 1$ for $e \in A$, and
$x(e) = 0$ otherwise.  At each step, we keep track of such a vector~$x_t$
and a set $E_t$ of edges, namely, the edges on which $x_t$ is random; all other
edges will be colored~$0$.  We choose an oriented edge $e = (v,w) \in E \setminus
E_t$, until such an edge $e$ no longer exists.  We set $x_{t+1}(e) = 1$ with
probability $p$, and $x_{t+1}(e) = 0$ with probability $1 - p$.  If $v$ and $w$
are already connected in~$x_t$ [i.e.,\ in the graph $(V, A_t)$ where $A_t
= \{e': x_t(e') = 1\} \subseteq E_t$], then we accept the edge and set $E_{t+1} =
E_t \cup \{e\}$. If $v$ and $w$ are not already connected, then we always accept
$x_{t+1}(e) = 0$, but we accept $x_{t+1}(e) = 1$ only with probability $1 / q$
(since by adding this edge we reduce by~$1$ the number of connected components).

When we reject, we know that $v$ and $w$ lie in separate components
in $(V, A_t)$.  To counteract this knowledge, to form $E_{t+1}$
we remove from $E_t$ all the
edges in the component of $(V, A_t)$ that contains~$w$, together with all
edges of~$E_t$ that lead out of this component (and which therefore do not
belong to~$A_t$).

We could cease our handling of a rejection step at this point and prove that
(a)~the algorithm works correctly and (b)~Theorem~\ref{THM:RC} below holds
(and the proof simplifies somewhat) with the bound on~$p$ decreased to
$$
p < 1 / (\Delta - (1/q)).
$$
However, we shall omit the formal proof of correctness and instead discuss a
small (provably valid) trick which gains us some efficiency.

Suppose that there are~$M$ vertices in the removed component.
Consider the (connected!)\ graph consisting of the vertices and edges in this
component, together with the vertex~$v$ and the edge~$\{v, w\}$.
Choose (in any fashion) a spanning tree~$T$ of this graph;
$T$ will comprise $M + 1$ vertices and therefore~$M$ edges.
Add back all these $M$ edges to get $E_{t+1}$.
Sample from the random cluster model on~$T$, and add back in the edges
thereby colored~$1$ to get $A_{t+1}$.

The key observation here is that it is elementary to sample from
the random cluster model when the graph is a tree.  Indeed, then
each edge independently is colored~$1$ with probability $\rho / (1 - p + \rho)$
and~$0$ with probability $(1 - p) / (1 - p + \rho)$, where
$$
\rho := p / q.
$$  

The random cluster model is an extension of the ferromagnetic Ising
and Potts models.
When $q > 1$ is an integer, and $p = 1 - \exp\{-\beta\}$, then samples from
the random cluster model may be used to generate samples from the
ferromagnetic Potts model with $q$ colors by independently taking each connected
component of $(V, A)$, uniformly choosing one of the~$q$ colors, and assigning to
every vertex in the component that color.   For certain instances of the random
cluster model, the heat bath Markov chain approach is believed from experimental
evidence to be rapidly mixing~\cite{proppw96}, but
no theoretical rapid mixing results in the positive direction are known for any
nontrivial instances of the problem.  For some instances, the Markov chain
approach is known \emph{not} to be rapidly mixing~\cite{gorej97}.
For the RR approach, we know that when $p$ is small
(corresponding to small $\beta$), the approach takes an
expected number of steps which is linear in the number of edges:

\begin{thm}
\label{THM:RC}
Suppose that
$$p < \frac{\Delta - (1/q) - \sqrt{[\Delta - (1/q)]^2 - 4[1 - (1/q)](\Delta - 1)}}
  {2 [1 - (1/q)] (\Delta - 1)}.$$
Then the expected number of steps required by the RR algorithm is $O(|E|)$.
\end{thm}

For example, if $\Delta = 4$ (as on a $2$-dimensional rectangular grid) and $q =
2$ (corresponding to the Ising model), then our restriction is that $p < 1/3$;
this improves on the restriction $p < 1 / (\Delta - (1/q)) = 2/7$ obtained when
the ``add a tree'' trick is not employed.

Comments analogous to those following the proof of Theorem~\ref{THM:WIS} again
apply.
\smallskip

\begin{proof}
We use a potential function that rewards us for adding edges and
penalizes us for connecting components.  Let
$$\phi(E_t, A_t) := |E_t| - \alpha c(A_t),$$
where $\alpha$ will be determined later.

When the edge $\{v, w\}$ we attempt to add to~$E_t$ is between two vertices
already connected in~$A_t$, then $\phi$ always goes up by 1, making this case
uninteresting. It is when $\{v, w\}$ would connect two previously unconnected
components of~$A_t$ that the calculation becomes interesting.

If the edge is chosen to be excluded from $A_{t+1}$, then $\phi$ increases by
$1$.  If the edge is proposed to be included in $A_{t+1}$, then $\phi$ changes by
$1 - \alpha$ if we accept.  If we reject, we remove from~$A_t$ (and also
from~$E_t$) a component of size $M$ and (from~$E_t$) all of its adjacent edges. 
Not counting the edge $\{v, w\}$ and making sure that we do not double-count, this
totals at most $M (\Delta - 1)$ edges removed from~$E_t$.  However, we add exactly
$M - 1$ new components to~$A_t$ by removing these edges.  When we add the tree~$T$
back in, this produces~$M$ new edges for~$E_{t+1}$, but for each such edge
there is a $\rho / (1 - p + \rho)$ chance of including the edge in $A_{t+1}$
and thereby reducing the number of components by~$1$.  Therefore, 
when we attempt to add $\{v, w\}$ to $A_{t+1}$, but reject instead, the 
expected contribution to the change in~$\phi$ is at least
$$-M(\Delta - 1) + M + \alpha\left(M - 1 - 
  M \left(\frac{\rho}{1 - p + \rho} \right) \right).
$$
Now $M$ may be very large (nearly as large as $n$), so we choose $\alpha$
in such a way that the coefficient of~$M$ in this expression vanishes.  That is,
we set
$$\alpha := (\Delta - 2) \left( \frac{1 - p + \rho}{1 - p} \right),$$
and so the contribution in this case is bounded below by $-\alpha$.

We try to put the edge in with probability $p$ and to leave it out with
probability $1 - p$.  We accept an inclusion with probability $1/q$.  
Putting everything together, 
we find that the expected change in $\phi$ at any time
step when $v$ and $w$ are not already connected in~$A_t$ is at least
$$(1 - p) 
 + p \left[\frac{1}{q}(1 - \alpha) + 
\left(1 - \frac{1}{q}\right)(-\alpha) \right],$$  
which is positive exactly when 
$$p <
  \frac{\Delta - (1/q) - \sqrt{[\Delta - (1/q)]^2 - 4[1 - (1/q)](\Delta - 1)}}
  {2 [1 - (1/q)] (\Delta - 1)}.$$
\end{proof}
$\Box$


In this case, the Markov chain approach does not have theoretical
guarantees on the running time for any nontrivial value of $p$.
While coupling from the past may also be used to generate 
perfect samples, there is no a priori bound on its running
time.

As with CFTP, we may still use the RR approach for values of $p$ for
which no theoretical bound exists.  We simply do not know beforehand
how long the algorithm will take.  Unlike CFTP, the RR approach is
interruptible, so we may abort the procedure if it needs too many steps,
without introducing bias into the sample. 

\paragraph{Proper colorings of a graph}
Finding the number of proper colorings of a graph is
a
\#P-complete problem~\cite{jerrum95}.  Recall that
a proper coloring of a graph assigns each vertex a color such that
no edge has both endpoints colored the same color.  The ability to
sample from the set of proper colorings leads to an approximation
algorithm for counting the number of such colorings.  

Markov chain approaches require that $k$, the number of colors,
be at least $(11/6)\Delta$
(where $\Delta$ is again the maximum 
degree of the graph) \cite{vigoda99} in order to guarantee rapid mixing
for the chain. 
Perfect sampling using bounding chains~\cite{huber98,haggstromn98} is
only guaranteed to run in polynomial time when the number of colors is
$\Omega(\Delta^2)$. Unfortunately, the straighforward RR approach does not match
these bounds.  Somewhat roughly stated,

\begin{thm}
The heat bath RR approach to generating perfect colorings requires only a linear
expected number of steps when $k$ is $\Omega(\Delta^4)$.
\end{thm}

As with the bounding chain procedure, however, this algorithm may
be run even when $k$ is much smaller; we simply have no reasonable a priori
bound on the running time in such cases.

\paragraph{The Move Ahead 1 chain}
Finally, we present a problem where an RR-based algorithm seems
experimentally to run
fast although we cannot give any theoretical bounds.
In the list update problem, a set of items is kept in a list.  To
access an item, a user starts at the beginning of the list
and steps through the items until the desired item is located.  The located item
may be replaced in the list anywhere between its current position and the front of
the list, at fixed cost.  The goal is to use a replacement strategy that keeps
small the access times (i.e., item depths in the list) needed for items.

Call the strategy which moves the accessed item to the front of the list
the \emph{Move to Front} (MTF) \emph{rule}.
A worst-case analysis shows that the MTF rule yields a 2-approximation for the
optimal total access time for any sequence of item requests~\cite{sleatort85}.  
Alternatively, it is useful to employ probabilistic models to describe how list
items are chosen to be accessed. Commonly, such an access model will induce a
Markov chain model on the evolution of the order of the list.  Characteristics
such as the limiting distribution as $t \to \infty$ of~$A_t$, where $A_t$ is the
access time for the item accessed at time~$t$, can then be estimated by drawing
from the stationary distribution of the chain.

To be specific, label the items with identification numbers $1, \ldots,
n$; suppose that at each time step, independently of previous time steps, any
particular item~$i$ is accessed with probability $p_i > 0$ (independently of the
order of the list); and suppose that after each selection, the accessed item is
moved forward one rank in the list, i.e.,\ is transposed with its predecessor in
the list.  (If the accessed item is already at the front of the list, the order
of the list is left unchanged.)  The self-organization rule we have described
is called the \emph{Move Ahead}~$1$ (MA1) \emph{rule}.  The limiting expected
access time for MA1 is known to be, for any access probability vector~$p$, no
more than that for MTF~\cite{rivest76}.  Further Monte Carlo study of the
limiting access time distribution is complicated by the fact that sampling from
the limiting list-order distribution~$\pi$ (for which a formula is known, but
only up to a normalizing constant) seems to be 
quite 
difficult in general.

Coupling from the past approaches to sampling from~$\pi$ exist~\cite{huber99}, but
experimental evidence suggests that use of RR gives a faster algorithm.
Suppose that 
$p_i \propto r^i$ for some ratio $0 < r \leq 1$.  Then experimental evidence
suggests that for each fixed value of $r \in (0, 1]$ the expected running time is
linear in~$n$, although the constant of linearity does varywith $r$.
The Markov chain approach to this problem is only known to be rapidly mixing when
$r < 0.2$~\cite{huber99}.


\section{Conclusion}
\label{SEC:Conc}
The RR approach to perfect sampling gives exact samples 
from difficult distributions without 
using 
the traditional
Markov chain.  It is quite different from other recent 
approaches to perfect sampling such as coupling from the past.

Because it dispenses with the Markov chain, the RR approach yields,
for restricted versions of some of these problems,
the first expected linear time algorithms for these problems.
Even when the running time of RR is unknown, the algorithm may be
run and the output will be guaranteed to come from the correct distribution. 

Unlike coupling from the past, RR is interruptible, so 
the user may set a time limit on
the algorithm's running time (if measured in number of iterations of
the basic Repeat loop)
without introducing bias into the
sample.  Like read-once coupling from the past~\cite{wilson2000},
this algorithm does not require storage of any random bits.
(Another perfect sampling approach, that of Fill, Machida, Murdoch,
and Rosenthal~\cite{fill99} is also interruptible but not read-once, and so does 
requires storage of random bits).  We wish to stress that these
existing means for perfect sampling rely on finding a ``good'' Markov
chain for the problem at hand.  RR does away with the chain, and in doing
so breaks the $O(n \ln n)$ barrier that has characterized so many of
these problems.

For independent sets and for proper colorings, the theoretical bounds
obtained apply only for a more restricted set of parameters than
do those based on Markov chain approaches.  However, when the 
appropriate restriction is met, our RR
method is faster, yielding samples in a linear (expected) number of steps.  
Moreover, much work has gone into analyses of Markov chains,
while our work is still rather new, and we might hope with time and further
effort eventually to match or even to relax the restrictions needed for the
Markov chain approaches.  For the Move Ahead 1 chain we do not know any
theoretical bounds on the running time of our method.   
However, 
computer experiments show that 
for this problem the RR method works much better in practice than
does the CFTP method.

For the random cluster model, our RR
technique is guaranteed to run in a linear (expected) number of steps for
a range of values of $p$.  This is in sharp constrast to the Markov chain
approach, where no polynomial running time bounds are known except in
trivial cases.

In summary, the randomness recycler is not applicable in all situations where
Markov chain approaches are used, but
RR often gives a fast read-once interruptible means
for generating perfect samples that in restricted cases gives the
first linear time algorithms for some difficult and important problems.

\bibliographystyle{plain}
\bibliography{c:/papers/huber/references/myrefs}

\end{document}